\documentclass[11pt,a4paper]{smfart} 
\usepackage{mathsmf}
\usepackage{smfthm}

\usepackage{graphicx}
\usepackage{color}
\usepackage[utf8]{inputenc}
\usepackage[T1]{fontenc}
\usepackage[french]{babel}
\usepackage{hyperref}

\title{Spectre et géométrie conforme 
des variétés compactes à bord}
\author{Pierre Jammes}
\address{Univ. Nice Sophia Antipolis, CNRS, LJAD, UMR 7351\\
06100 Nice France}
\email{pjammes@unice.fr}
\date{}
\begin{document}
\begin{abstract}
On montre que sur toute variété compacte $M^n$ à bord, il
existe une classe conforme $C$ telle que pour toute métrique $g\in C$
de volume~1, la première valeur propre non nulle du laplacien de Neumann 
vérifie $\lambda_1(M^n,g)< n\Vol(S^n,g_{\textrm{can}})^{2/n}$. On montre
une majoration analogue pour la première valeur propre non nulle 
de Steklov. La démonstration est basée sur une 
décomposition en anse de la variété. On montre aussi que le volume conforme
de $(M,C)$, l'invariant de Friedlander-Nadirashvili et le volume
de Möbius de $M$ sont ceux de la sphère $(S^n,g_{\textrm{can}})$.
Si $M$ est un domaine euclidien, hyperbolique ou sphérique, alors $C$ 
est la classe conforme de la métrique canonique. 
\end{abstract}
\begin{altabstract}
We prove that on any compact manifold $M^n$ with boundary, there exist a 
conformal class $C$ such that for any riemannian metric $g\in C$ of
unit volume, the first positive eigenvalue of the Neumann Laplacian satisfies 
$\lambda_1(M^n,g)< n\Vol(S^n,g_{\textrm{can}})^{2/n}$. We also prove a 
similar inequality for the first positive Steklov eigenvalue.
 The proof 
relies on a handle decomposition of the manifold. We also 
prove that the conformal
volume of $(M,C)$ is $\Vol(S^n,g_{\textrm{can}})$, and that the 
Friedlander-Nadirashvili and the Möbius volume of $M$ are equal
to those of the sphere.  If $M$ is a domain in a space form, $C$ is 
the conformal class of the canonical metric.
\end{altabstract}

\altkeywords{first Neumann and Steklov eigenvalue, conformal volume.}

\maketitle

\section{Introduction}
Considérons une variété riemannienne compacte $(M^n,g)$, avec ou sans bord,
et notons $\lambda_1(M,g)>0$ la première valeur propre non nulle
du laplacien sur $M$, en posant la condition de Neumann s'il y a un bord.
Il découle des travaux de J.~Hersch \cite{he70} que sur la sphère de
dimension~2, on a
\begin{equation}
\lambda_1(S^2,g)\Vol(S^2,g)\leq8\pi,
\end{equation}
l'égalité ayant lieu pour la métrique canonique.
Mais on sait que cette inégalité ne s'étend pas aux autres surfaces closes
(les premiers exemples ont été donnés par P.~Buser \cite{bu84}. Voir
aussi \cite{bm01} et les références qui y sont données)
On sait que la borne supérieure de $\lambda_1(M^2,g)\Vol(M^2,g)$ vaut
$12\pi$ sur le plan projectif (\cite{ly82}), $8\pi^2/\sqrt3$ sur le tore
(\cite{na96}, \cite{gi09}), que sur la bouteille de Klein elle s'exprime
en fonction d'une intégrale elliptique (\cite{jnp06}, \cite{esgj06}),
et on conjecture qu'elle vaut $16\pi$ sur la surface orientable de genre~2 
(\cite{jlnnp05}).
On peut en déduire facilement
à l'aide de \cite{ces03} que cette borne supérieure est minorée par $12\pi$
sur les surfaces closes autres que la sphère.
En dimension supérieure ou égale à~3,
on sait que $\lambda_1(M^n,g)\Vol(M^n,g)^{2/n}$ n'est pas borné
(voir \cite{ta79}, \cite{mu80}, \cite{cd94}).

 Récemment, G.~Kokarev et N.~Nadirashvili ont donné une généralisation
de l'inégalité de Hersch sur les surfaces orientables à bord sous des
hypothèses conformes. Si $g$ est une métrique sur $M$, la classe
conforme de $g$ est définie par $[g]=\{hg, h\in C^\infty, h>0\}$.
\begin{theo}[\cite{kn10}]\label{intro:kn}
Soit $M^2$ une surface orientable compacte à bord. Il existe une
classe conforme $C$ sur $M^2$ telle que pour toute métrique $g\in C$,
on a $\lambda_1(M^2,g)\Vol(M^2,g)\leq8\pi$.
\end{theo}

 Le but de cet article est de démontrer une généralisation de ce résultat
aux variétés compactes à bord (orientables ou non) de dimension quelconque.
On va aussi montrer un énoncé du même type pour la première valeur propre
du spectre de Steklov. Rappelons que si $\rho\in C^\infty(\partial M)$
est une fonction positive sur le bord de $M$, le spectre de Steklov est la
suite $0=\sigma_0(M,g,\rho)<\sigma_1(M,g,\rho)\leq\sigma_2(M,g,\rho)\ldots$
des réels $\sigma$ pour lesquels le problème de Steklov, défini par les
équations
$\Delta f=0$ dans $M$ et $\frac{\partial f}{\partial \nu}=\sigma \rho f$
sur $\partial M$, admet une solution non nulle (voir le
paragraphe~\ref{conforme:stek} pour plus de détails).

\begin{theo}\label{intro:th1}
Soit $M^n$ une variété compacte à bord de dimension~$n\geq2$. Il existe
une classe conforme $C$ sur $M^n$ telle que pour toute métrique $g\in C$,
\begin{equation}\label{intro:eq1}
\lambda_1(M,g)\Vol(M,g)^{2/n}< n\omega_n^{2/n}
\end{equation}
et
\begin{equation}\label{intro:eq2}
\sigma_1(M,g,\rho)\mathcal M(\partial M)\Vol(M)^{\frac{2-n}n}<
n\omega_n^{2/n},
\end{equation}
où $\omega_n=\Vol(S^n,g_{\textrm{can}})$ et
$\mathcal M(\partial M)=\int_{\partial M}\rho\ \de v_g$ désigne
la masse totale de la densité $\rho$.

Si $M$ est un domaine euclidien, sphérique ou hyperbolique,  alors ces
inégalités sont vraies pour la classe conforme de la métrique canonique.
\end{theo}
\begin{rema}
L'inégalité~(\ref{intro:eq1}) est optimale : on peut
assez facilement faire tendre la première valeur propre vers celle
de la sphère à volume fixé dans une classe conforme donnée quelconque,
et donc $\sup_{g\in C}\lambda_1(M,g)\Vol(M,g)^{2/n}\geq n\omega_n^{2/n}$ en
général
(la démonstation donnée dans~\cite{ces03} pour les variétés closes
s'applique aussi aux variétés à bord).
La situation est différente pour l'inégalité~(\ref{intro:eq2}). Par exemple,
selon l'inégalité de Weinstock \cite{we54}, qui est un analogue de
l'inégalité de Hersch pour le problème de Steklov, on a
$\sigma_1(M,g,\rho)\mathcal M(\partial M)\leq2\pi$ si $M$ est homéomorphe
à un disque.
\end{rema}
\begin{rema}\label{intro:rkpe}
Le théorème~\ref{intro:th1} contraste fortement avec les résultats
obtenus récemment par R.~Petrides sur les variétés closes. En effet,
il montre dans \cite{pe14} et \cite{pe13} que pour toute classe conforme
$C$ sur une variété close $M^n$, on a
$\sup_{g\in C}\lambda_1(M,g)\Vol(M,g)^{2/n}>n\omega_n^{2/n}$, sauf si
$(M^n,C)$ est conformément équivalente à la sphère ronde.
\end{rema}
\begin{rema}
 En dimension~2, l'inégalité~(\ref{intro:eq1}) améliore celle de
Kokarev et Nadirashvili sur deux points : elle est stricte et on n'a pas
besoin de supposer que la surface est orientable.
\end{rema}
\begin{rema}
Si $M$ est un domaine de $\R^n$, $S^n$ ou $H^n$, on obtient une majoration
de $\lambda_1(M,g)$ en fonction du volume qui est plus faible que celles
de Szegö-Weinberger ou de Ashbaugh-Benguria \cite{ab95} mais qui vaut pour
toute métrique conforme. Dans ce contexte, on peut interpréter
l'inégalité~(\ref{intro:eq1}) comme une inégalité de Szegö-Weinberger conforme.
Par ailleurs, Colbois, El Soufi et Girouard ont donné dans \cite{cesg11}
une majoration de $\sigma_k(M)\mathcal M(\partial M)
\Vol(M)^{\frac{2-n}n}$ sur les domaines de $\R^n$, $S^n$ ou $H^n$ qui ne
dépend que de $n$ et $k$. La majoration~(\ref{intro:eq2}) ne vaut que
pour $k=1$, mais elle a l'avantage d'être plus explicite et
conformément invariante.
\end{rema}

L'inégalité~(\ref{intro:eq1}) est spécifique à la condition
de Neumann. Avec la condition de Dirichlet, on peut faire tendre
la première valeur propre vers l'infini dans une classe conforme à volume fixé:
\begin{theo}\label{intro:th3}
Soit $M^n$ une variété compacte à bord de dimension~$n\geq 2$. Pour toute
classe conforme $C$ sur $M^n$, on a
$$\sup_{g\in C}\lambda_1^D(M^n,g)\Vol(M^n,g)^{2/n}=+\infty$$
où $\lambda_1^D(M^n,g)$ désigne la première valeur propre du laplacien
avec condition de Dirichlet.
\end{theo}
On montrera en fait ce résultat pour le laplacien agissant sur les
formes différentielles de degré $p\geq2$, le laplacien de Neumann
correspondant au cas $p=0$ et le laplacien de Dirichlet au cas $p=n$.

 Comme l'article \cite{kn10} utilise des techniques spécifiques à la
dimension~2 (invariance conforme de la norme $L^2$ du gradient, fonctions
méromorphes sur des courbes hyperelliptiques), la généralisation du
théorème~\ref{intro:kn} en grande dimension peut surprendre. Le principe de la
démonstration, déjà utilisé dans \cite{ja08}, est de faire appel à la
théorie du volume conforme développée dans \cite{ly82} et \cite{esi86}
et d'utiliser une décomposition en anses de la variété. Dans le cas
du spectre de Steklov, on passera par l'intermédiaire du volume
conforme relatif introduit par A.~Fraser et R.Schoen dans \cite{fs11}.

 La même technique va nous permettre de calculer certains invariants des
variétés compactes à bord. Rappelons d'abord quelques définitions. Le
volume conforme est un invariant conforme des variétés compactes
défini par
\begin{equation}\label{intro:vc}
V_c(M,[g])=\inf_{\stackrel{m\geq 1}{\varphi:(M,[g])\hookrightarrow S^m}}
\sup_{\gamma\in G_m}\Vol(\gamma\circ\varphi(M)),
\end{equation}
l'application $\varphi$ parcourant l'ensemble des immersions
conformes de $(M,[g])$ dans $S^m$, et $G_m$ désignant le groupe de Möbius de
dimension~$m$, c'est-à-dire le groupe des difféomorphismes conformes de $S^m$.
Nous utiliserons le fait que le volume conforme permet de majorer
la première valeur propre du laplacien (\cite{ly82}, \cite{esi86}):
\begin{equation}\label{intro:vc2}
\lambda_1(M,g)\Vol(M,g)^{2/n}\leq n V_c(M,[g])^{2/n}.
\end{equation}

Dans \cite{fn99}, L.~Friedlander et N.~Nadirashvili ont défini un
nouvel invariant différentiable des variétés compactes en posant
\begin{equation}\label{intro:nu}
\nu(M)=\inf_g\sup_{\tilde g\in[g]}\lambda_1(M,\tilde g)
\Vol(M,\tilde g)^{2/n}.
\end{equation}
Enfin, j'ai défini dans \cite{ja08} le volume de Möbius d'une variété compacte
par
\begin{equation}\label{intro:vm}
V_\mathcal{M}(M)=\inf_{\stackrel{m}{\varphi:M\hookrightarrow S^m}}
\sup_{\gamma\in G_m}\Vol(\gamma\circ\varphi(M)),
\end{equation}
la différence avec le volume conforme étant que $\varphi$ parcourt
l'ensemble de toutes les immersions de $M$ dans $S^m$ sans restriction.
Il découle de l'inégalité (\ref{intro:vc2}) que $\nu(M)\leq
nV_\mathcal{M}(M)^{2/n}$. On peut alors montrer :

\begin{theo}\label{intro:th2}
Soit $M^n$ une variété compacte à bord. Alors
\begin{enumerate}
\item Il existe une classe conforme $C$ sur $M$ telle que $V_c(M,C)=\omega_n$;
\item $\nu(M)=n\omega_n^{2/n}$;
\item $V_\mathcal{M}(M)=\omega_n$.
\end{enumerate}
De plus, si $M$ est un domaine euclidien, sphérique ou hyperbolique, alors
$C$ est la classe conforme de la métrique canonique.
\end{theo}
Ces résultats contrastent avec le cas des variétés compactes sans bord.
En effet, pour les variétés closes, on a toujours $V_c(M,C)\geq\omega_n$
et il découle des travaux de R.~Petrides  que l'égalité caractérise la sphère
conformément ronde (cf. remarque~\ref{intro:rkpe}).
On sait qu'on a aussi $\nu(M)\geq n\omega_n^{2/n}$ mais il existe très
peu de variétés sur lesquelles on sait calculer cet invariant, les
seuls exemples sont $\nu(S^n)=n\omega_n^{2/n}$, $\nu(\R P^2)=12\pi$
(cela découle du fait que $\R P^2$ n'admet qu'une seule classe conforme),
 $\nu(T^2)=\nu(K^2)=8\pi$ où $K^2$ désigne la bouteille de Klein
(voir~\cite{gi09}) et les surfaces compactes orientables en général, pour
lesquelles l'invariant $\nu$ vaut encore $8\pi$
(\cite{pe14}). Enfin, la seule variété close dont on connaît le
volume de Möbius est la sphère. En général, on sait seulement que
$V_\mathcal{M}$ est minoré par $\omega_n$ et uniformément majoré à dimension
fixée (voir~\cite{ja08}).

Comme exemple d'application du théorème~\ref{intro:th2},
on peut déduire une majoration de la deuxième valeur propre des
opérateurs de Schrödinger sur les domaines euclidiens, sphériques et
hyperboliques en utilisant le résultat d'El~Soufi et Ilias \cite{esi92}:

\begin{cor}\label{intro:cor1}
Soit $M$ une variété conforme à un domaine compact de $\R^n$ (par exemple
domaine de $S^n$ ou $H^n$) et $V\in C^\infty(M)$.
La deuxième valeur propre de l'opérateur de Schrödinger $\Delta+V$
sur $M$ avec condition de Neumann est majorée par
$n(\frac{\omega_n}{\Vol(M)})^{2/n}+\overline V$,
où $\overline V$ désigne la valeur moyenne de $V$.
\end{cor}

 On commencera par quelques rappels sur le volume conforme et ses
liens avec le spectre de Neumann et de Steklov dans la
section~\ref{conforme}, on démontrera au passage les
théorèmes~\ref{intro:th1} et~\ref{intro:th2} dans le cas des domaines
de $\R^n$, $S^n$ et $H^n$.
On traitera ensuite le cas des surfaces,
pour lequel les arguments topologiques sont plus simples qu'en dimension
quelconque, dans la
section~\ref{2}. Puis nous verrons le cas général dans la section~\ref{n}.
Enfin, la dernière section sera consacrée à la démonstration
du théorème~\ref{intro:th3}.

Je remercie B. Colbois d'avoir porté l'article \cite{kn10} à ma connaissance,
ainsi que le rapporteur de l'article, dont les remarques ont permis d'améliorer
le texte.

\section{Spectre et volume conforme}%
\label{conforme}
\subsection{Première valeur propre de Neumann}
Les démonstrations des théorèmes~\ref{intro:th1} et~\ref{intro:th2}
reposent sur les propriétés des immersions conformes dans la sphère
et du volume conforme étudiées dans \cite{ly82} et \cite{esi86}.
On utilisera en particulier le résultat suivant, en notant
$G_m$ le groupe de Möbius de dimension~$m$, c'est-à-dire
le groupe des difféomorphismes conformes de la sphère~$S^m$. 
\begin{lemme}[\cite{esi86}]\label{conf:lemme}
Soit $(M^n,g)$ une variété compacte et $\varphi:M\mapsto S^m$ une immersion
conforme de $M$ dans la sphère de dimension~$m$. Il existe un élément $\gamma$
du groupe de Möbius $G_m$ tel que $\lambda_1(M,g)\Vol(M,g)^{2/n}\leq
n\Vol(\gamma\circ\varphi(M))^{2/n}$.
\end{lemme}
\begin{rema}\label{conf:rem1}
Ce lemme n'est pas énoncé explicitement dans \cite{ly82} ou \cite{esi86},
mais sa démonstration est contenue dans celle du théorème~2.2 de \cite{esi86}.
Par ailleurs, cette majoration est obtenue en appliquant le principe
du min-max à des fonctions test qui ne vérifient \emph{a priori} aucune
condition de bord. Elle est donc valable pour le laplacien avec 
la condition de Neumann mais pas la condition de Dirichlet.
\end{rema}
\begin{rema}\label{conf:rem2}
Dans \cite{ly82} et \cite{esi86}, le volume de $\gamma\circ\varphi(M)$
est majoré par sa borne supérieure sous l'action de $G_m$. Pour les
immersions que nous allons considérer dans cet article, cette borne 
supérieure, qui vaudra $\omega_n$, ne sera pas atteinte. C'est la raison 
pour laquelle les inégalités du théorème~\ref{intro:th1} sont strictes.
\end{rema}

On est ainsi ramené au problème de construire une immersion de $M\to S^m$ 
dont le volume reste strictement inférieur à $\omega_n$ sous l'action 
du groupe de Möbius.

Pour illustrer l'efficacité de cette technique, nous allons montrer
dès à présent le théorème~\ref{intro:th1} pour  la première valeur
propre de Neumann (on verra comment adapter la démonstration au cas de la 
première valeur propre de Steklov au paragraphe suivant, cf. 
remarque~\ref{conf:stekrem}),
ainsi que le théorème~\ref{intro:th2}, dans 
le cas d'un domaine de l'espace euclidien, hyperbolique ou de la sphère:
\begin{theo}\label{conf:th}
Soit $M$ un domaine (strict) de $\R^n$, $S^n$ ou $H^n$. Alors
\begin{enumerate}
\item la classe conforme $C$ de la métrique canonique sur $M$ vérifie 
$V_c(M,C)=\omega_n$ et pour toute métrique $g\in C$,
$\lambda_1(M^n,g)\Vol(M^n,g)^{2/n}< n\omega_n^{2/n}$;
\item $\nu(M)=n\omega_n^{2/n}$;
\item $V_\mathcal{M}(M)=\omega_n$.
\end{enumerate}
\end{theo}
\begin{proof}
Dans le cas d'un domaine de $S^n$, l'immersion considérée est simplement
l'identité. Sous l'action du groupe de Möbius, le domaine $M$ reste
un domaine strict, donc son volume reste strictement inférieur à $\omega_n$.
On déduit du lemme~\ref{conf:lemme} que $\lambda_1(M,g)\Vol(M,g)^{2/n}<
n\omega_n^{2/n}$ pour toute métrique $g\in C$. Par définition du volume 
conforme, on en déduit aussi que $V_c(M,C)\leq\omega_n$. Les autres résultats
du théorème découlent des inégalités générales $\omega_n\leq V_{\mathcal M}(M)
\leq V_c(M,C)$ et $n\omega_n^{2/n}\leq\nu(M)\leq n V_{\mathcal M}(M)^{2/n}$
montrées dans~\cite{esi86}, \cite{fn99} et~\cite{ja08}.

Si $M$ est un domaine de $\R^n$ la projection stéréographique 
donne une immersion conforme dans $S^n$, et si $M$ est un domaine
de $H^n$ on utilise le fait que $H^n$ est conforme à un hémisphère
de $S^n$. Le reste de la démonstration est identique.
\end{proof}

\subsection{Première valeur propre de Steklov}\label{conforme:stek}
Le problème des valeurs propres de Steklov consiste à résoudre l'équation
\begin{equation}
\left\{\begin{array}{ll}
\Delta f=0 & \textrm{dans }M\\
\frac{\partial f}{\partial \nu}=\sigma \rho f & \textrm{sur }\partial M
\end{array}\right.
\end{equation}
où $\nu$ est un vecteur unitaire sortant normal au bord et 
$\rho\in C^\infty(\partial M)$ une fonction densité fixée. On note
$\mathcal M(\partial M)=\int_{\partial M}\rho\de v_g$ (si $\rho=1$ on parle
de problème de Steklov homogène, et on a dans ce cas $\mathcal M(\partial M)
=\Vol(\partial M,g)$). L'ensemble des
réels $\sigma$ solutions du problème forme un spectre discret positif
noté
\begin{equation}
0=\sigma_0(M,g,\rho)<\sigma_1(M,g,\rho)\leq\sigma_2(M,g,\rho)\ldots.
\end{equation}
Dans le cas homogène, le spectre de Steklov est aussi connu comme étant 
le spectre de l'opérateur Dirichlet-to-Neumann. Le problème de Steklov,
déjà étudié à la fin du XIX\ieme{} siècle et au début du XX\ieme{} 
(voir~\cite{st99}, \cite{st02} et les références qui y sont données), 
apparaît dans divers problèmes physiques. Par exemple il
permet de modéliser l'évolution d'une membrane libre dont la masse se
concentre sur son bord, et il intervient dans certains problèmes de 
tomographie. Concernant les bornes conformes des $\sigma_i$ pour tout $i$,
on peut consulter~\cite{ha11}.

Dans~\cite{fs11}, A.~Fraser et R.~Schoen ont montré la majoration conforme
suivante, dans le cas homogène:
\begin{equation}\label{conf:stek}
\sigma_1(M,g)\Vol(\partial M,g)\Vol(M)^{\frac{2-n}n}\leq nV_{rc}(M)^{2/n}
\end{equation}
où $V_{rc}(M)$ est un invariant conforme, baptisé volume conforme relatif, 
défini par
\begin{equation}\label{conf:vrc}
V_{rc}(M,[g])=\inf_{\stackrel{m\geq 1}{\varphi:(M,[g])\hookrightarrow B^{m+1}}}
\sup_{\gamma\in G_m}\Vol(\gamma\circ\varphi(M)),
\end{equation}
où $B^{m+1}$ est la boule unité de $\R^{m+1}$ et $\varphi$ parcourt l'ensemble 
des immersions conformes dans la boule $B^{m+1}$ telles que 
$\varphi(\partial M)\subset S^m$. On utilise le fait
que l'action du groupe de Möbius sur la sphère $S^m$ s'étend naturellement
en une action conforme sur $B^{m+1}$. 

En considérant le cas particulier des
immersions $\varphi:M\to S^m$, il est clair que $V_{rc}(M)\leq V_c(M)$.
On peut donc déduire immédiatement une majoration de $\sigma_1(M)$ 
en fonction du volume conforme dans le cas homogène. On va maintenant
montrer l'inégalité~(\ref{conf:stek}) dans le cas non homogène:
\begin{theo}
Soit $(M,g)$ une variété compacte à bord et $\rho\in C^\infty(\partial M)$
une fonction. Si on note $\sigma_1(M,g,\rho)$ la première valeur propre
non nulle du problème de Steklov avec densité $\rho$, alors
$$\sigma_1(M,g,\rho)\mathcal M(\partial M)\Vol(M)^{\frac{2-n}n}\leq 
nV_{rc}(M)^{2/n}.$$
\end{theo}
\begin{proof}
Rappelons d'abord que la valeur propre $\sigma_1(M,g,\rho)$ possède la 
caractérisation variationnelle suivante (cf. \cite{ba80}) :
\begin{equation}\label{conf:minmax}
\sigma_1(M,g,\rho)=\inf_{f\neq0}\frac{\int_M|\de f|^2\de v_g}%
{\int_{\partial M}f^2\rho\ \de v_g}\ \textrm{ où }
\int_{\partial M}f\rho\ \de v_g=0.
\end{equation}

Soit $\varphi : M\to B^{m+1}$ une immersion conforme telle que
$\varphi(\partial M)\subset S^m$. Selon le lemme~2.1 de \cite{esi92}, 
il existe un élément $\gamma\in G_m$ tel que si on note $\tilde\varphi=
\gamma\circ\varphi$, les coordonnées $\tilde\varphi_i$ de  $\tilde\varphi$
vérifient $\int_{\partial M}\tilde\varphi_i\rho\ \de v_g=0$.
En appliquant le principe variationnel~(\ref{conf:minmax}) aux fonctions 
$\tilde\varphi_i$ on obtient $\sigma_1(M,g,\rho)
\int_{\partial M}\tilde\varphi_i^2\rho\ \de v_g\leq
\int_M|\de\tilde\varphi_i|^2\de v_g$. Après sommation, en utilisant le fait
que $\sum_i\tilde\varphi_i^2=1$ sur $\partial M$ et en appliquant
une inégalité de Hölder, on obtient par un calcul devenu classique
\begin{equation}\label{conf:stekeq1}
\sigma_1(M,g,\rho)\leq\frac{\sum_i\int_M|\de\tilde\varphi_i|^2\de v_g}%
{\sum_i\int_{\partial M}\tilde\varphi_i^2\rho\ \de v_g}\leq
\frac{\left(\int_M\left(\sum_i|\de\tilde\varphi|^2\right)^{n/2}\de v_g\right)
^{2/n}}{\mathcal M(\partial M)\Vol(M,g)^{(2-n)/n}}.
\end{equation}
Comme l'immersion $\tilde\varphi$ est conforme, on a aussi 
$\tilde\varphi^*(g_\textrm{eucl})=
\frac1n\left(\sum_i|\de\tilde\varphi|^2\right)g$ et donc 
\begin{equation}\label{conf:stekeq2}
\left(\int_M\left(\sum_i|\de\tilde\varphi|^2\right)^{n/2}\de v_g\right)^{2/n}
=n\Vol(\tilde\varphi(M))^{2/n}=n\Vol(\gamma\circ\varphi(M))^{2/n}.
\end{equation}
Par définition, $\Vol(\gamma\circ\varphi(M))$ est majoré par $V_{rc}(M)$,
ce qui permet de conclure.
\end{proof}
\begin{rema}\label{conf:stekrem}
L'inégalité~(\ref{intro:eq2}) du théorème~\ref{intro:th1}
se démontre de la même manière que l'inégalité~(\ref{intro:eq1}) en remplaçant
le lemme~\ref{conf:lemme} par les inégalités~(\ref{conf:stekeq1})
et~(\ref{conf:stekeq2}). En particulier, le cas des domaines euclidiens
se traite comme dans le théorème~\ref{conf:th}.
\end{rema}

\section{Première valeur propre des surfaces à bord}\label{2}

Comme dans \cite{ja08}, nous allons commencer par reformuler le problème
en projetant stéréographiquement la sphère $S^m$ sur $\R^m$.
Rappelons d'abord que si $\gamma$ est un élément quelconque de $G_m$ et
$r$ une isométrie de la sphère $S^m$, alors $\Vol(r\circ\gamma\circ\varphi
(M))=\Vol(\gamma\circ\varphi(M))$ pour toute immersion $\varphi$ de $M$
dans $S^m$. On peut donc se ramener à étudier l'action d'un sous-ensemble
$G'_m$ de $G_m$ (pas nécessairement un
sous-groupe) tel que tout élément de $G_m$ s'écrive sous la forme
$r\circ\gamma$ avec $\gamma\in G'_m$ et $r\in\SO(m)$.

Il existe plusieurs choix possibles pour $G'_m$, mais l'un d'entre
eux est particulièrement adapté à la démonstration qui suit. En projetant
stéréographiquement la sphère $S^m$ sur $\R^m\cup\{\infty\}$, on se
ramène à considérer des immersions $\varphi:M\hookrightarrow
\R^m\cup\{\infty\}$ en calculant les volumes pour la métrique
$g_{S^m}=\frac4{(1+\|x\|^2)^2}g_{\mathrm{eucl}}$ où $g_{\mathrm{eucl}}$
désigne la métrique euclidienne canonique. On peut alors choisir
pour $G'_m$ l'ensemble des homothéties et des translations de
$\R^m$ (voir le théorème~3.5.1 de~\cite{be83}).

La difficulté par rapport au cas des domaines traités dans la section précédente
est de trouver une immersion de $M$ dans une sphère telle que le volume
de l'immersion
reste strictement inférieure à $\Vol(S^n,g_{\textrm{can}})$ sous
l'action du groupe de Möbius. Le reste de la démonstration est identique.

\begin{figure}[h]
\begin{center}
\begin{picture}(0,0)%
\includegraphics{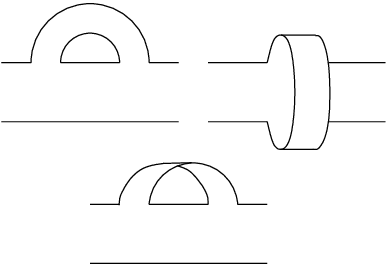}%
\end{picture}%
\setlength{\unitlength}{4144sp}%
\begingroup\makeatletter\ifx\SetFigFont\undefined%
\gdef\SetFigFont#1#2#3#4#5{%
  \reset@font\fontsize{#1}{#2pt}%
  \fontfamily{#3}\fontseries{#4}\fontshape{#5}%
  \selectfont}%
\fi\endgroup%
\begin{picture}(2949,2001)(664,-1603)
\end{picture}%
\end{center}
\caption{Adjonctions d'anses à une surface à bord\label{2:2anse}}
\end{figure}

On va donc montrer l'existence d'une telle immersion dans le cas où $n=2$
en utilisant une décomposition en anses de la surface à bord. Précisons
qu'une adjonction d'anse à une surface à bord consiste à coller un
produit $[0,1]\times[0,1]$ sur la surface en identifiant $\{0\}\times[0,1]$
et $\{1\}\times[0,1]$ avec deux intervalles du bord de la surface. Les
adjonctions d'anses sont de trois types topologiques différents,
représentés sur la figure~\ref{2:2anse}, selon que le nombre
de composantes connexes du bord augmente de 1, diminue de 1 ou reste constant
(auquel cas la surface obtenue est non orientable).

La décomposition en anse peut alors s'énoncer ainsi:
\begin{lemme}
Toute surface compacte à bord de caractéristique d'Euler $\chi$ s'obtient
topologiquement à partir du disque par adjonction de $1-\chi$ anses.
\end{lemme}
\begin{proof}
On se base sur le fait qu'une surface compacte est entièrement
déterminée par sa caractéristique d'Euler, le nombre de composantes
du bord et son orientabilité (\cite{hi94}, théorème~3.1). Remarquons aussi
qu'ajouter une anse diminue de~1 la caractéristique d'Euler.

Commençons par le cas des surfaces orientables. Pour obtenir une surface
de caractéristique $\chi$ et ayant $p$ composantes de bord, il suffit
d'attacher au disque $k=(p-\chi)/2$ anses qui créent une composante
de bord, puis $k'=1-(p+\chi)/2$ anses
qui connectent deux composantes de bord différentes (pour une surface
orientable, $p$ et $\chi$ ont la même parité, donc $k$ et $k'$
sont bien entiers). La caractéristique
d'Euler de la surface obtenue est bien $1-k-k'=\chi$ et le nombre de
composantes connexes du bord est $1+k-k'=p$.

Pour obtenir une surface non orientable de caractéristique $\chi$
et ayant $p$ composantes de bord, on attache $p-1$ anses qui créent une
composante de bord, puis $2-\chi-p$ anses qui ne changent pas le nombre
de composantes du bord.
\end{proof}

Nous allons maintenant montrer l'existence d'une immersion dont
l'aire est contrôlée sous l'action du groupe de Möbius.

\begin{theo}\label{2:th}
Soit $M$ une surface compacte à bord. Il existe une
immersion $\varphi:M\to S^3$ telle
que $\Vol(\gamma\circ\varphi(M))<4\pi$ pour tout
élément $\gamma$ du groupe de Möbius~$G_3$ de~$S^3$.
\end{theo}

\begin{proof}
Étant donnée une décomposition en anse de la surface, et deux paramètres
$\delta,\varepsilon>0$ petits, on se donne un plongement de la surface
dans $\R^3$ tel que les anses aient une largeur (pour la métrique
euclidienne $g_{\mathrm{eucl}}$) égale à $\varepsilon$ et que la courbure
du plongement (plus précisément la norme de la seconde forme fondamentale)
soit majorée  par $\delta$ (les valeurs de $\delta$ et
$\varepsilon$ seront fixées plus loin). On lisse
le bord à la jonction des anses et on note
$\varphi:M\hookrightarrow \R^3$ le plongement obtenu. La
figure~\ref{2:surface} représente un tel plongement pour un tore privé de
deux disques.
\begin{figure}[h]
\begin{center}
\begin{picture}(0,0)%
\includegraphics{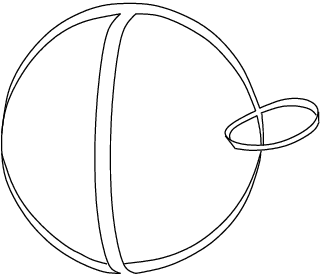}%
\end{picture}%
\setlength{\unitlength}{4144sp}%
\begingroup\makeatletter\ifx\SetFigFont\undefined%
\gdef\SetFigFont#1#2#3#4#5{%
  \reset@font\fontsize{#1}{#2pt}%
  \fontfamily{#3}\fontseries{#4}\fontshape{#5}%
  \selectfont}%
\fi\endgroup%
\begin{picture}(2441,2081)(1568,-2355)
\put(3236,-1116){\makebox(0,0)[lb]{\smash{{\SetFigFont{11}{13.2}{\rmdefault}{\mddefault}{\updefault}{\color[rgb]{0,0,0}(a)}%
}}}}
\put(2446,-1591){\makebox(0,0)[lb]{\smash{{\SetFigFont{11}{13.2}{\rmdefault}{\mddefault}{\updefault}{\color[rgb]{0,0,0}(b)}%
}}}}
\end{picture}%
\end{center}
\caption{Surface à anses fines\label{2:surface}}
\end{figure}
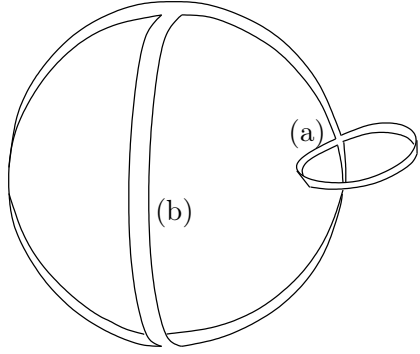

La clef de la démonstration est que si $\gamma\in G'_3$,
l'aire des parties de $\gamma\circ\varphi(M)$ projetées près de $\infty$
devient négligeable pour la métrique $g_{S^3}$, même si le rapport
d'homothétie de $\gamma$ est grand (cf.~\cite{ja08}).
En effet, une homothétie de rapport $R$ va multiplier
l'aire euclidienne par $R^2$, mais dans la formule
$g_{S^3}=\frac4{(1+\|x\|^2)^2}g_{\mathrm{eucl}}$ la norme de $x$ est de
l'ordre de $R$ pour les parties de $\varphi(M)$
projetées près de l'infini et donc leur aire devient petite quand $R$
devient grand (de l'ordre de $R^{-2}$). On va maintenant préciser cette
majoration.

Quand $\varepsilon$
est petit, l'aire de $\gamma\circ\varphi(M)$ est  négligeable sauf si
le rapport d'homothétie de $\gamma$ est grand et qu'une partie
de $\gamma\circ\varphi(M)$ reste près de l'origine. On va donc se placer
dans le cas où $\gamma$ est une homothétie centrée en un point $p$ de
la surface et dont le rapport $R$ est grand, et étudier l'aire de
la surface quand $R\to\infty$.

On considère la boule $B$ de centre $p$ et de rayon~1 et on sépare
la surface en deux morceaux $\varphi(M)\cap B$ et $\varphi(M)\backslash B$.
Pour $\delta$ et $\varepsilon$ suffisamment petits, il existe une
constante $c_1>0$ telle que l'aire
$V_1=\Vol(\gamma\cdot(\varphi(M)\backslash B))$ vérifie $V_1\leq c_1R^{-2}$
pour $R$ assez grand. En outre, pour $R$ fixé, l'aire $V_1$ est
asymptotiquement proportionnelle à $\varepsilon$ quand $\varepsilon$ tend
vers~0.
La constante $c_1$ tend donc vers~0 quand $\varepsilon\to0$.

 Pour estimer l'aire $V_2=\Vol(\gamma\cdot(\varphi(M)\cap B))$,
on considère la projection orthogonale $\pi$ de la surface sur son plan
tangent $\mathcal P$ en $p$. Cette projection prend par exemple l'une des
formes représentées sur la figure~\ref{2:zoom} selon que $p$ est près d'une
jonction d'anse ou pas. Le cas où $p$ est situé sur le bord de la surface
sera précisé à la fin de la démonstration.

\begin{figure}[h]
\begin{center}
\begin{picture}(0,0)%
\includegraphics{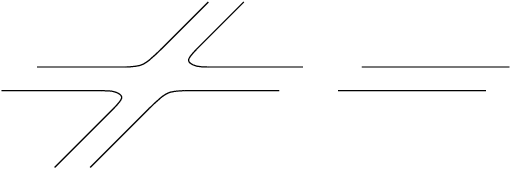}%
\end{picture}%
\setlength{\unitlength}{4144sp}%
\begingroup\makeatletter\ifx\SetFigFont\undefined%
\gdef\SetFigFont#1#2#3#4#5{%
  \reset@font\fontsize{#1}{#2pt}%
  \fontfamily{#3}\fontseries{#4}\fontshape{#5}%
  \selectfont}%
\fi\endgroup%
\begin{picture}(3894,1284)(1249,-1558)
\put(2566,-1186){\makebox(0,0)[lb]{\smash{{\SetFigFont{12}{14.4}{\rmdefault}{\mddefault}{\updefault}{\color[rgb]{0,0,0}(a)}%
}}}}
\put(4186,-1186){\makebox(0,0)[lb]{\smash{{\SetFigFont{12}{14.4}{\rmdefault}{\mddefault}{\updefault}{\color[rgb]{0,0,0}(b)}%
}}}}
\end{picture}%
\end{center}
\caption{Voisinage d'un point de la surface\label{2:zoom}}
\end{figure}
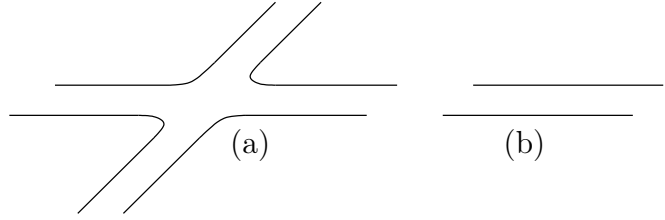

Le paramètre $\delta$ majore la courbure de $\varphi(M)$, donc celle
de $\gamma\cdot\varphi(M)$ est majorée par $\delta/R$ (pour la métrique
euclidienne).
Si $\delta$  est suffisamment petit, l'angle entre le plan tangent en un point
$x\in
\gamma\cdot(\varphi(M)\cap B)$ proche de $p$ et le plan $\mathcal P$ est
majoré (au premier ordre par rapport à $\|x-p\|$) par $\frac\delta R\|x-p\|$.
Le jacobien de $\pi$
au point $x$ sera donc encadré grossièrement par $1\pm
2\frac\delta R\|x-p\|$. Cette estimation est faite pour la métrique
euclidienne, mais la densité de la mesure sphérique décroit à la même
vitesse en $x$ et $\pi(x)$ quand $x$ s'éloigne de $p$, donc elle reste
valable pour la métrique sphérique.

On a donc une majoration de la forme
\begin{equation}
|V_2'-V_2|\leq
2\int_{\pi\circ\gamma\cdot(\varphi(M)\cap B)}\frac\delta R\|x-p\|\de x
\end{equation}
où $V_2'=\Vol(\pi\circ\gamma\cdot(\varphi(M)\cap B))$, l'intégrale étant
calculée pour la métrique sphérique. Comme le rapport
$\|x-p\|/R$ reste borné, on peut trouver
une constante $c_2$ telle que l'intégrale du membre de droite soit
majorée (pour la métrique sphérique) par
$c_2\delta R^{-2}$ pour $R$ assez grand. En outre, pour $R$ fixé, la largeur
du domaine d'intégration $\pi\circ\gamma\cdot(\varphi(M)\cap B)$
est proportionnelle à $\varepsilon$, donc la constante $c_2$ tend vers~0
quand $\varepsilon$ tend vers~0.

Pour conclure, on va distinguer le cas où $p$ est dans l'intérieur de
la surface et le cas où $p$ est sur le bord. Considérons d'abord
le cas où $p$ est à l'intérieur. La projection stéréographique envoie le
plan $\mathcal P$ sur une sphère ronde plongée dans $S^3$, son aire
pour la métrique sphérique est donc majorée par celle d'une sphère
équatoriale, c'est-à-dire $4\pi$. Par conséquent l'aire $V_2'$ est majorée
par $4\pi-\Vol(U)$,
où $U$ est le complémentaire de $\pi\circ\gamma\cdot(\varphi(M)\cap B)$
dans $\mathcal P$. L'aire $\Vol(U)$ est de l'ordre de $R^{-2}$ quand
$R\to\infty$, et elle croît quand $\varepsilon$ décroît.

Finalement, l'aire $V_1+V_2$ est majorée par $c_1R^{-2}+c_2\delta R^{-2}+
4\pi-\Vol(U)$. Comme les constantes $c_1$ et $c_2$ tendent vers~0 quand
$\varepsilon\to0$, on peut bien trouver un $\varepsilon$ tel que
$c_1R^{-2}+c_2\delta R^{-2}-\Vol(U)$ soit toujours strictement négatif.
On obtient que $\Vol(\gamma\circ\varphi(M))<4\pi$ pour tout $R$.

Quand $p$ est sur le bord, la projection de la surface sur $\mathcal P$
tend vers un demi-plan quand $R\to\infty$. Par conséquent, le domaine
$U$ défini précédemment tend vers le demi-plan complémentaire dans
$\mathcal P$, et en particulier son aire ne tend pas vers~0. La majoration
précédente est donc \emph{a fortiori} vraie dans ce cas.

Les estimations précédentes dépendent \emph{a priori} du point $p$,
mais par compacité de la surface le choix des constantes peut être fait
uniformément sur la surface. Finalement, on a bien
$\Vol(\gamma\circ\varphi(M))<4\pi$ pour tout $\gamma$.
\end{proof}

\section{Variétés à bord de dimension quelconque}\label{n}
La démonstration en dimension quelconque utilise des arguments géométriques
semblables à ceux de dimension~2, mais les arguments topologiques sont plus
techniques. Dans le même esprit que~\cite{ja08}, on fera appel à une
décomposition en anses de la variété. On se référera aux livres \cite{ra02}
(chapitres~1, 2 et 6) et \cite{ko07} (chapitre~VI et~VII) pour les
aspects topologiques de la démonstration.

Soit $M$ une variété compacte à bord et $S^{k-1}\hookrightarrow\partial M$
une sphère plongée dans le bord de $M$ dont le fibré normal (dans
$\partial M$) est trivial. Un voisinage tubulaire de $S^{k-1}$ dans
$\partial M$ est alors un produit trivial $S^{k-1}\times B^{n-k}$.
Ajouter une $k$-anse le long de $S^{k-1}$ consiste à coller une variété
produit $B^k\times B^{n-k}$ (dont le bord est la réunion des deux variétés
$B^k\times S^{n-k-1}$ et $S^{k-1}\times B^{n-k}$ recollées le long
de $S^{k-1}\times S^{n-k-1}$)  sur un voisinage tubulaire
de la sphère $S^{k-1}$. C'est un résultat classique de la théorie du
cobordisme qu'on peut obtenir n'importe quelle variété connexe compacte
à bord par adjonctions successives d'un nombre fini d'anses à partir d'une
boule (\cite{ra02} théorème~2.2 ; \cite{ko07} ch.~VII, théorème~1.1 et~1.2).
En outre, on peut
faire en sorte qu'on n'ait pas à
utiliser de $n$-anse, c'est-à-dire à \og{}~boucher~\fg{} un bord en y collant
une boule (\cite{ko07}, ch.~VII, théorème~6.1). La démonstration des
théorèmes~\ref{intro:th1} et~\ref{intro:th2} va
consister à raisonner par récurrence sur le nombre d'anses.
Comme en dimension~2, il suffit de montrer l'existence d'une immersion
dont le volume est uniformément contrôlé sous l'action du groupe de
Möbius:
\begin{theo}\label{n:th}
Soit $M$ une variété compacte à bord de dimension $n$. Il existe une
immersion $\varphi:M\to S^m$ pour un entier $m$ suffisamment grand telle
que $\Vol(\gamma\circ\varphi(M))<\Vol(S^n,g_{\textrm{can}})$ pour tout
élément $\gamma$ du groupe de Möbius~$G_m$ de~$S^m$.
\end{theo}
\begin{proof}
Si $M$ est une boule, n'importe quel plongement  $\varphi:M\to S^n$
convient.

Supposons qu'une variété à bord $M$ vérifie la conclusion du théorème ;
on va montrer qu'elle reste vraie pour une variété $\tilde M$ obtenue
en ajoutant une $k$-anse qu'on notera $\mathcal A$.

On note $\varphi$ l'immersion de $M$ dans $S^m$ telle que
$\Vol(\gamma\circ\varphi(M))<\Vol(S^n,g_{\textrm{can}})$ pour tout $\gamma$.
Comme dans le paragraphe précédent, on projette stéréographiquement
$S^m$ sur $\R^m\cup\infty$ et on se restreint aux éléments $\gamma\in G'_m$.
On va construire une immersion $\tilde\varphi$ de $\tilde M$ en collant une
anse
$\mathcal A=B^k\times B^{n-k}$ sur $\varphi(M)$ (pour ne pas alourdir
le texte, on identifiera les variétés à leur immersion).
Le long de la sphère d'attache $S^{k-1} \hookrightarrow\partial M$, on
commence par coller une boule $B^k$ de manière à ce que le long de
$S^{k-1}$, l'espace tangent à $B^k$ soit contenu dans l'espace tangent de
$M$. Quitte à augmenter la dimension de la sphère $S^m$,
on peut faire en sorte que $B^k$ ne rencontre pas $M$ (c'est-à-dire
que $M\cup B^k$ est plongée). L'anse $B^k\times B^{n-k}$ est ensuite
obtenue en épaississant $B^k$ et en lissant $\tilde M$ le long de l'attache
de l'anse (voir figure~\ref{n:anse}). Si le rayon $\varepsilon$
(pour la métrique
euclidienne de $R^m$) de $B^{n-k}$ est suffisamment petit, l'immersion
$\tilde\varphi$ de $\tilde M$ ainsi obtenue est un plongement.
Comme dans le cas des surfaces, on choisira plus loin un $\varepsilon$
négligeable devant la courbure de la boule $B^k$.

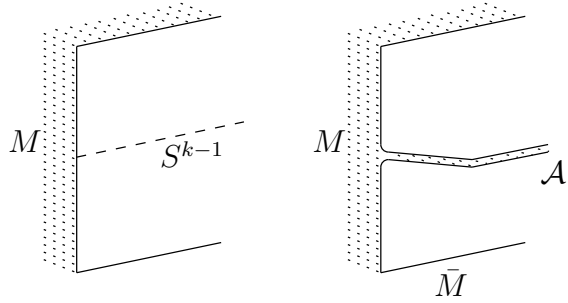
\begin{figure}[h]
\begin{center}
\begin{picture}(0,0)%
\includegraphics{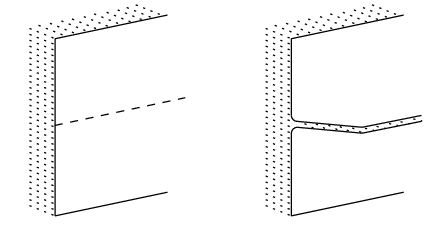}%
\end{picture}%
\setlength{\unitlength}{4144sp}%
\begingroup\makeatletter\ifx\SetFigFont\undefined%
\gdef\SetFigFont#1#2#3#4#5{%
  \reset@font\fontsize{#1}{#2pt}%
  \fontfamily{#3}\fontseries{#4}\fontshape{#5}%
  \selectfont}%
\fi\endgroup%
\begin{picture}(3228,1831)(481,-1385)
\put(1396,-556){\makebox(0,0)[lb]{\smash{{\SetFigFont{12}{14.4}{\rmdefault}{\mddefault}{\updefault}{\color[rgb]{0,0,0}$S^{k-1}$}%
}}}}
\put(2296,-466){\makebox(0,0)[lb]{\smash{{\SetFigFont{12}{14.4}{\rmdefault}{\mddefault}{\updefault}{\color[rgb]{0,0,0}$M$}%
}}}}
\put(496,-466){\makebox(0,0)[lb]{\smash{{\SetFigFont{12}{14.4}{\rmdefault}{\mddefault}{\updefault}{\color[rgb]{0,0,0}$M$}%
}}}}
\put(3016,-1321){\makebox(0,0)[lb]{\smash{{\SetFigFont{12}{14.4}{\rmdefault}{\mddefault}{\updefault}{\color[rgb]{0,0,0}$\bar M$}%
}}}}
\put(3646,-646){\makebox(0,0)[lb]{\smash{{\SetFigFont{12}{14.4}{\rmdefault}{\mddefault}{\updefault}{\color[rgb]{0,0,0}$\mathcal A$}%
}}}}
\end{picture}%
\end{center}
\caption{Voisinage de l'attache d'une anse\label{n:anse}}
\end{figure}

Il reste à choisir $\varepsilon$ de manière à contrôler le volume de ce
plongement quand on fait agir le groupe $G'_m$. Il est clair que sous
l'action des translations le volume de $M$ atteint un maximum (en effet, le
volume de $\gamma\cdot M$ tend alors vers zéro à l'infini pour la
métrique $g_{S^m}$). Ce maximum est
strictement plus petit que $\Vol(S^n,g_{\textrm{can}})$ par hypothèse,
et le volume de l'anse reste majoré par une borne qui tend vers~0
quand $\varepsilon$ tend vers~0 (il est essentiel ici pour contrôler
le volume de la $k$-anse que $k$ soit différent de $n$). On peut donc
choisir $\varepsilon$
de sorte que le volume de $\tilde M$ reste strictement plus petit
que $\Vol(S^n,g_{\textrm{can}})$ sous l'action des translations.

On est donc ramené à l'étude du volume sous l'action des homothéties.
On considère une homothétie $\gamma$ centrée en un point $p$ et de rapport $R$.
Le problème est de majorer le volume de $\tilde M$ quand $p\in\tilde M$
et $R$ est grand (si $R$ est inférieur à une borne donnée, on peut rendre
le volume de l'anse négligeable et conclure comme dans le cas des
translations).
On va en fait montrer par récurrence un résultat plus précis, à savoir
qu'il existe une constante $c_M$ telle que
\begin{equation}\label{n:majR}
\Vol(\gamma\cdot M)\leq\Vol(S^n,g_{\textrm{can}})-c_MR^{-n}
\end{equation}
quand $R\to\infty$.
On va distinguer quatre cas, selon que $p$ est sur $M$,
dans l'anse $\mathcal A$, sur le bord de l'anse ou près de la jonction
de $\mathcal A$ et $M$.

Supposons que $p$ soit dans $M$. Par hypothèse de récurrence,
$M$ vérifie la majoration~(\ref{n:majR}). Quand $R\to\infty$, l'anse
$\mathcal A$ tend vers l'infini et son volume pour la métrique
sphérique est majorée par $cR^{-n}$ où $c$ est une constante qui tend
vers~0 quand $\varepsilon\to0$. On a donc $\Vol(\gamma\cdot \bar M)\leq
\Vol(\gamma\cdot M)+cR^{-n}\leq\Vol(S^n,g_{\textrm{can}})-(c_M-c)R^{-n}$.
Pour $\varepsilon$ suffisamment petit, la constante $(c_M-c)$ est
positive ce qui permet de conclure. On peut noter que cette démonstration
reste valable si $p$ est sur le bord de $M$.

Si $p$ est dans l'intérieur de l'anse $\mathcal A$, on raisonne comme en
dimension~2.
On note $\mathcal A'$ un petit voisinage de $p$
on note $\pi$ la projection orthogonale sur l'espace tangent à
$\mathcal A$ en $p$, et par le même calcul que dans la section précédente
on obtient l'existence
d'une constante $c$ (dépendant de la courbure de l'anse) telle que
\begin{equation}
|\Vol(\pi(\gamma\cdot \mathcal A'))-\Vol(\gamma\cdot \mathcal A')|\leq cR^{-n}.
\end{equation}
Comme $\gamma\cdot M$ et $\gamma\cdot (\mathcal A\backslash\mathcal A')$
tendent vers le point à l'infini quand $R\to\infty$,
il existe une constante $c$ telle que $\Vol(\gamma\cdot M)+\Vol(\gamma\cdot
(\mathcal A\backslash\mathcal A'))\leq cR^{-n}$.
En notant $U$ le complémentaire de $\pi(\gamma\cdot \mathcal A')$ dans
le plan tangent on a $\Vol(\pi(\gamma\cdot \mathcal A'))=
\Vol(S^n,g_{\textrm{can}})-\Vol(U)$, avec $\Vol(U)\sim c'R^{-n}$ quand
$R\to\infty$ pour une constante $c'$ qu'on peut rendre arbitrairement grande
en faisant tendre $\varepsilon$ vers~0 (c'est-à-dire qu'on peut choisir
à partir de quelle valeur de $R$ ce volume devient négligeable). On a
finalement $\Vol(\gamma\cdot\bar M)\leq
\Vol(S^n,g_{\textrm{can}})-c''R^{-n}$ avec une constante $c''$ qui est bien
positive si $\varepsilon$ est suffisamment petit.

Si $p$ est sur le bord de l'anse $\mathcal A$, on modifie la démonstration
précédente comme en dimension~2 : $\gamma\cdot\bar M$ tend vers un
demi-espace de $T_p\bar M$ quand $R$ tend vers l'infini, donc le volume de
$U$ ne tend pas vers~0 et la conclusion reste vraie.

Si $p$ est près de la jonction avec l'anse $\mathcal A$,
la démonstration est la même que si $p\in\mathcal A$.

 Comme en dimension~2,
on peut choisir le paramètre $\varepsilon$ uniforme par rapport à $p$ par
compacité.
\end{proof}

\section{Grandes valeurs propres du laplacien de Dirichlet}
On va maintenant démontrer le théorème~\ref{intro:th3}. On va en fait
prouver un résultat plus général, à savoir qu'on peut faire 
tendre toutes les valeurs propres du laplacien de Hodge-de~Rham vers
l'infini en fixant le volume et la classe conforme, à l'exception
des valeurs propres du laplacien de Neumann.

Rappelons quelques définitions. Le laplacien de Hodge-de~Rham, qui 
agit sur les formes différentielles de la variété $M$, est défini par
\begin{equation} 
\Delta=\de\codiff+\codiff\de:\Omega^p(M)\to\Omega^p(M),
\end{equation}
où $\de$ désigne la différentielle extérieure et $\codiff$ sont 
adjoint $L^2$.  On considérera la condition de bord suivante, dite absolue, 
en notant $j$ l'inclusion $\partial M\hookrightarrow M$ et $N$ un champ
de vecteur normal au bord: 
\begin{equation}
(A)\ \left\{\begin{array}{l}j^*(\prodint_N\omega)=0\\
j^*(\prodint_N\de\omega)=0.\end{array}\right.
\end{equation}

On notera $0=\lambda_{p,0}(M,g)<\lambda_{p,1}(M,g)
\leq\lambda_{p,2}(M,g)\leq\ldots$ son spectre en restriction aux $p$-formes,
en répétant les valeurs propres non nulles s'il y a multiplicité 
(la multiplicité de la valeur propre nulle, si elle existe, est un 
invariant topologique: c'est le nombre de Betti $b_p(M)$).

Avec la condition de bord (A), le laplacien de Hodge-de~Rham restreint
aux $0$-formes, c'est-à-dire aux fonctions, est le laplacien de Neumann.
En degré $p=n$, les valeurs propres sont celles du laplacien de Dirichlet.
En effet, si $f$ est une fonction et $\de v_g$ la forme volume, alors
$\Delta(f\de v_g)=(\Delta f)\de v_g$, et $f\de v_g$ vérifie la condition
(A) si et seulement si $f$ vérifie la condition de Dirichlet.

\begin{theo}
Soit $M^n$ une variété compacte à bord de dimension~$n\geq 2$. Pour toute
classe conforme $C$ sur $M^n$, il existe une famille de métriques $g_t\in C$,
$t>0$ telle que, pour tout $p\geq2$, on a 
$\lambda_{p,1}(M,g_t)\Vol(M,g_t)^{2/n}\to+\infty$ quand $t$ tend vers 
$+\infty$.
\end{theo}
\begin{rema}
On ne peut pas faire tendre $\lambda_{1,1}(M,g_i)\Vol(M,g_i)^{2/n}$ 
vers l'infini. En effet, comme la différentielle extérieure commute
avec le laplacien, le spectre restreint aux 1-formes contient le 
spectre du laplacien de Neumann.
\end{rema}
\begin{rema}
Sur les variétés closes, on sait qu'on peut faire tendre 
$\lambda_{p,1}(M)\Vol(M)^{2/n}$ vers l'infini en fixant le volume et 
la classe conforme, mais seulement pour $p\in[2,n-2]$ (voir~\cite{ces06}).
La méthode, devenue classique pour faire tendre les valeurs propres
de certains opérateurs vers l'infini, consiste à déformer une petite 
boule de la variété en un cylindre très long (cf.~\cite{gp95}, \cite{ab00},
 \cite{ces06}, \cite{aj11}). Nous allons ici
adapter cette technique à la présence d'un bord.
\end{rema}
\begin{proof}
On procédera en trois étapes. D'abord, on se ramènera au cas
où la métrique est euclidienne au voisinage d'un point du bord. 
On construira ensuite une famille de métriques qui fait tendre le volume
vers l'infini et on montrera enfin que le spectre de la variété est
uniformément minoré pour cette famille de métriques.

La première étape se base sur le résultat de J.~Dodziuk \cite{do82} selon 
lequel, si deux métriques sont proches pour la distance de Lipschitz, 
alors leurs spectres sont proches aussi. Plus précisément, si deux métriques 
$g$ et $\tilde g$ vérifient $\tau^{-1}g\leq\tilde g\leq\tau g$, 
alors $\tau^{-(3n-1)}\lambda_{p,1}(M,\tilde g)\leq\lambda_{p,1}(M,g)\leq
\tau^{3n-1}\lambda_{p,1}(M,\tilde g)$ et $\tau^{-n/2}\Vol(M,\tilde g)\leq
\Vol(M,g)\leq\tau^{n/2}\Vol(M,\tilde g)$. Si $h\in C^\infty(M)$ est 
une fonction strictement positive et que $\tau^{-1}g\leq\tilde g\leq\tau g$,
alors on a aussi $\tau^{-1}hg\leq\tilde hg\leq\tau hg$. On peut en déduire
que si on peut faire tendre $\lambda_{p,1}(M,hg)\Vol(M,hg)^{2/n}$
vers l'infini en faisant varier $h$, alors $\lambda_{p,1}(M,h\tilde g)
\Vol(Mh\tilde g)^{2/n}$ tend aussi vers l'infini. Autrement dit, il suffit
de montrer le théorème~\ref{intro:th3} pour une classe conforme particulière,
le résultat général s'en déduit.

On peut donc se placer dans le cas où la classe conforme contient
une métrique $g$ qui est euclidienne dans un voisinage d'un point du bord.
On peut alors adapter la construction de~\cite{ces06}  dans ce contexte.
On sait qu'une boule euclidienne peut être déformée de manière conforme 
en la réunion d'un cylindre et d'un hémisphère. Plus
précisément, si on note $r$ la coordonnée radiale sur la boule euclidienne
$B(\varepsilon,g_{\textrm{eucl}})$ de rayon~$\varepsilon$ et qu'on pose
\begin{equation}
h_\varepsilon(r)=\left\{\begin{array}{ll}
\frac\varepsilon r&\textrm{si }\varepsilon e^{-\frac L\varepsilon}\leq
r\leq\varepsilon,\\
e^{\frac L\varepsilon}&\textrm{si }0\leq r\leq
\varepsilon e^{-\frac L\varepsilon},\\
\end{array}\right.
\end{equation}
alors la boule $B(\varepsilon,h_\varepsilon^2g_{\textrm{eucl}})$ est
isométrique à la réunion d'un cylindre de rayon~$\varepsilon$ et de 
longueur~$L$ et d'un hémisphère de rayon~$\varepsilon$.
Appliquée à une demi-boule centrée en un point du bord et en prolongeant 
la fonction $h_\varepsilon$ par~1 en dehors de cette demi-boule, cette 
déformation conforme crée un demi-cylindre de longueur~$L$ au bout duquel 
est collé un quart de sphère, comme sur la figure~\ref{dirichlet:nez}.
On note $g_L$ la métrique obtenue.
\begin{figure}[h]
\begin{center}
\begin{picture}(0,0)%
\includegraphics{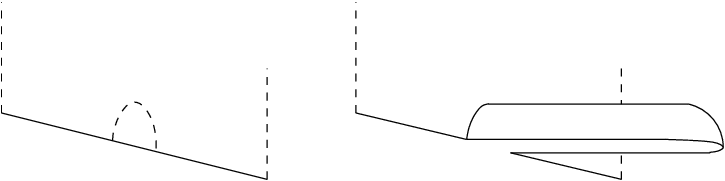}%
\end{picture}%
\setlength{\unitlength}{4144sp}%
\begingroup\makeatletter\ifx\SetFigFont\undefined%
\gdef\SetFigFont#1#2#3#4#5{%
  \reset@font\fontsize{#1}{#2pt}%
  \fontfamily{#3}\fontseries{#4}\fontshape{#5}%
  \selectfont}%
\fi\endgroup%
\begin{picture}(5523,1374)(664,-748)
\end{picture}%
\end{center}
\caption{déformation de la métrique sur le bord de la variété%
\label{dirichlet:nez}}
\end{figure}

 On peut alors minorer le spectre de la variété obtenue par la même
technique que dans~\cite{gp95} et~\cite{ces06}. On recouvre la variété
$M$ par trois ouverts $U_i$, $i=1,2,3$ définis de la manière suivante :
$U_1$ est formé de la réunion de $M\backslash B(\varepsilon)$ et
d'une portion du demi-cylindre de longueur~1, l'ouvert $U_2$ est
formé du demi-cylindre de longueur $L$ et $U_3$ est formé du quart
de sphère et de d'une portion de demi-cylindre adjacente de longueur~1.

Les intersections $U_1\cap U_2$ et $U_2\cap U_3$ sont topologiquement 
le produit d'un intervalle et d'une boule, leur cohomologie est donc
triviale en degré $p\geq1$ (y compris en degré~$n$ et $n-1$, à la 
différence de la construction de~\cite{ces06}), ce qui permet d'appliquer 
le lemme de McGowan \cite{mc93} (plus précisément, on utilise l'énoncé
donné dans~\cite{gp95}, lemme~1) pour tous les degrés $p\geq2$ : 
il existe des constantes $a,b,c>0$ ne dépendant pas de $L$ telles que 
\begin{equation}
\lambda_{p,1}(M,g_L)\geq\frac a{\frac b{\mu_p(U_2)}+c}\textrm{ pour }p\geq2,
\end{equation}
où $\mu(U_2)$ désigne la plus petite valeur propre non nulle du laplacien
de Hodge sur $U_2$ en degré $p$ (avec la condition de bord absolue).

Comme l'ouvert $U_2$ est un produit riemannien, son spectre est déterminé
par la formule de Künneth et on peut montrer que $\mu(U_2)$ est minoré 
indépendamment de $L$ (le calcul est le même que dans~\cite{gp95}
et~\cite{ces06}).

On obtient finalement que $\lambda_{p,1}(M,g_L)$ reste uniformément
minoré pour $p\geq2$ quand $L\to+\infty$ tandis que $\Vol(M,g_L)\to
+\infty$, ce qui permet de conclure.
\end{proof}

\newcommand{\etalchar}[1]{$^{#1}$}
\providecommand{\bysame}{\leavevmode ---\ }
\providecommand{\og}{``}
\providecommand{\fg}{''}
\providecommand{\smfandname}{\&}
\providecommand{\smfedsname}{\'eds.}
\providecommand{\smfedname}{\'ed.}
\providecommand{\smfmastersthesisname}{M\'emoire}
\providecommand{\smfphdthesisname}{Th\`ese}

\end{document}